\documentclass[11pt,reqno,draft]{amsart}
\usepackage{latexsym,amssymb}

\newtheorem{defi}{Definition}[section]
\newtheorem{proposition}[defi]{Proposition}

\newtheorem{rem}[defi]{Remark}
\newtheorem{lemma}[defi]{Lemma}
\newtheorem{theorem}[defi]{Theorem}

\newcommand{\defeq}{\mathrel{\mathrm{\raise0.1ex\hbox{:}\hbox{=}\strut}}}
 
\DeclareMathOperator{\arctg}{arctg} \DeclareMathOperator{\s}{span}

\def\dd{\:\mathrm{d}}

\def\R{\mathbb R}
\def\NN{\mathbb{N}}

\def\E{\mathbb E}

\def\PP{\mathbb P}

\def\L{ L}

\def\la{\langle}
\def\ra{\rangle}

\def\ein{\mathbf{1}}
\def\Lip{\mathop{\rm Lip}}

\def\ub{\bar{u}}

\def\ve{\varepsilon}
\def\smallsum{\mathop{\mbox{$\sum$}}}
\def\px{\partial _x}

\def\qt{Q^T}
\def\dotsk{{i=1, \dots, k}}

\begin{document}

\numberwithin{equation}{section}

\title[Ergodicity of  Stochastic  Curve Shortening Flow in the Plane]{Ergodicity of  Stochastic  Curve Shortening Flow in the Plane}
%\title[Semilinear SPDE's with multiplicative noise ]{Solution for quasi linear stochastic parabolic problem with linear growth functionals  }
\author[A. Es--Sarhir]{Abdelhadi Es--Sarhir}
\author[M-K. von Renesse]{Max-K. von Renesse}
\address{Technische Universit\"at Berlin, Institut f\"ur Mathematik \newline Stra{\ss}e des 17. Juni 136, D-10623 Berlin, Germany}

\email{[essarhir,mrenesse]@math.tu-berlin.de}
\thanks{The authors acknowledge support
from the DFG Forschergruppe 718 "Analysis and Stochastics in Complex
Physical Systems".}

\def\krylovrozovskii{Krylov-Rozovski{\u\i}}
\def\subvstern{_{V^*}}

\keywords{Degenerate stochastic equations, monotonicity, \textit{e}-property.}


\begin{abstract}
We study  a model  of the motion by mean curvature
of an (1+1) dimensional interface in a 2D Brownian velocity field. For the well-posedness of the model  we
prove existence and uniqueness  for certain degenerate nonlinear
stochastic evolution equations in the variational framework of  \krylovrozovskii ,  replacing  the standard coercivity
assumption  by a Lyapunov type condition. Ergodicity is established for the  case of  additive noise,  using the lower bound technique for Markov  semigroups by Komorowski, Peszat and Szarek \cite{KPS}.
\end{abstract}

\maketitle
\section{Introduction}
\noindent   Motion by mean curvature is a well studied and rich
object in geometric PDE theory for which a variety of methods have
been developed (see  e.g.\ \cite{MR1931534} for a survey). In
physics it arises   as sharp interface limit of the Allen-Cahn
equation for the phase field of a binary alloy, describing the
motion of the interface between the two phases. Stochastic mean
curvature flow was derived heuristically  in e.g.,
\cite{citeulike:2163102} as a refined model incorporating the
influence of thermal noise. In the (d+1)-dimensional graph case the
corresponding SPDE is of the form
\begin{equation}
 d u = \sqrt {1 + |\nabla u|^2}\, {\rm div}\bigl( \frac{\nabla
u}{\sqrt {1 + |\nabla u|^2}} \bigr)\, dt+ B(u, \nabla u)  \delta W,  \label{grafspde}
\end{equation}
where $\delta$ stands for Stratonovich or It\^{o} differential,
depending on the model. The degeneracy of the drift operator makes  a rigorous treatment of this family of models very difficult.  Motivated by the deterministic theory Lions and Souganidis  introduced a notion  of  stochastic viscosity solutions    \cite{MR1659958,MR1799099},  but some technical details of this approach are still awaiting  full justification  \ \cite{MR1920103,caruana-2009}. Existence of weak subsequential limits  along  tight approximations of stochastic mean curvature flow has been obtained by Yip \cite{MR1656479} and more recently by R\"oger  and Weber \cite{webrog}.\\

In this paper we consider the  special case of a  (1+1)-dimensional   graph  interface in an e.g.\ 2D  Brownian velocity field, corresponding to the equation
 \begin{equation}   d u =   \frac{ \partial _x^2  u }{ 1 + (\partial _x u)^2} \, dt+ \sum\limits_{i=1}^{\infty} \phi_i(.,u(.))\,db_t^i. \label{themodel}
\end{equation}
\noindent In the deterministic case this equation is  also known as
curve shortening flow. Note that the  mild solution approach by da
Prato-Zabzcyk \cite{MR1207136} is not applicable because equation
\eqref{themodel} is not semilinear, i.e.\ does not contain  a
dominating linear component. For the analysis of \eqref{themodel}
we first establish an abstract existence and uniqueness result  in
the classical variational SPDE framework of  \krylovrozovskii{}
\cite{krylovrozovskii} for a certain class of nonlinear stochastic
evolution equations,  which are not coercive  but satisfy an
alternative  Lyapunov condition.  This is then applied to  equation
\eqref{themodel}  which is treated  in   the Gelfand triple
\[  H_0^{1} ([0,1]) \subset L^2([0,1]) \subset H^{-1}([0,1]),\]
although   the  operator $A: H_0^{1}([0,1])\to H^{-1}([0,1])$
\[ A u =  \frac{  \partial _x^2 u }{ 1 + (\partial _x u )^2} \]
fails to be coercive. By our method we  prove  well-posedness of
\eqref{themodel}, assuming  $u_0 \in H^{1}_0$, $\phi_i \in
\Lip([0,1]\times \R)$, $\phi_i(0,.)=\phi_i(1,.)=0$ and,  for some
finite $\Lambda$,
\begin{equation}
\sum_{i=1}^ \infty (\Lip(\phi_i))^2 \leq \Lambda ^2. \label{regcond} \end{equation}
The latter condition should be compared to the weaker assumption that for all $z_1, z_2 \in [0,1]\times \R$
 \begin{equation}
\sum_{i=1}^ \infty (\phi_i(z_1)-\phi_i(z_2))^2 \leq \Lambda^2 |z_1-z_2|^2, \label{kunitaregcond}
\end{equation}
which is well-known e.g.\ in the theory of isotropic flows, where it guarantees the existence of a forward stochastic flow $d\Phi = F(\Phi, dt)$ of homeomorphisms of $[0,1]\times \R_+$ driven by  the martingale field $F(z,t) = \sum _{i=1}^\infty \phi_i (z) b^{i}_t $, cf. \cite[Theorem 4.5.1]{MR1070361}.

\smallskip
In fact,  we show below that the SPDE \eqref{themodel} with noise
field satisfying only \eqref{kunitaregcond} and initial condition
$u_0 \in L^2([0,1])$ still admits a unique generalized solution
which is defined by approximation. More precisely, we obtain a
unique  Markov process $(\hat u_t^x; x\in L^2([0,1]); t\geq 0)$ on
$L^2([0,1])$, inducing a Feller semigroup on the space of bounded
continuos functions on $L^2([0,1])$ as the
unique generalized solution of \eqref{themodel}. However, in view of
the poor regularity of the operator $A$,  a more explicit
characterization of the $L^2([0,1])$-valued process $(\hat u
_t^x)_{t\geq0}$ by some SPDE or even just an associated Kolmogorov
operator on smooth finitely based test functions does not seem to be
available. This is very similar to the  generalized solutions for
abstract SPDE with only $m$-accretive drift operators obtained in
\cite{MR1400370} by   means of nonlinear semigroup theory. The
advantage in the present case is, however, that the variational
approach is embedded such that we know the solution $(\hat
u^x_t)_{t\geq 0}$ is strong if \eqref{regcond} holds and the initial
condition $x=u_0$ belongs to $H_0^1([0,1])$.
\smallskip

Finally we show the ergodicity of the generalized solution $(\hat u^x_t) $ of \eqref{themodel}  in the case of
 additive noise by  verifying the  conditions of  a recent abstract result by Komorowski, Peszat and Szarek for Markov semigroups
with the so-called \textit{e}-property \cite[Theorem 1]{KPS}.  We point out   that \cite[Theorem 3]{KPS} does not apply in our sitation because the deterministic flow does not converge to equilibrium locally uniformly with respect to the initial condition. However, for the verification of the lower bound in our case we exploit the fact that the stochastic flow admits a Lyapunov function with compact sublevel sets.

\section{Well-Posedness of certain non-coercive variational SPDE }
\subsection{\small \sc Strong solutions for a class of non-coercive SPDE with regular initial condition}
\noindent Although we are mainly interested in the example
\eqref{themodel} we shall formulate here a general existence and
uniqueness result in the abstract  variational framework of
\cite{krylovrozovskii} for stochastic evolution equations, following
with only a few changes  the excellent presentation in  \cite{Ro}.
Let \[ V \subset H\] be a continuous and dense embedding of two
separable Hilbert spaces with corresponding inner products $\langle
.,.\rangle_V$ and $\la .,.\ra_H$.  Via the Riesz isomorphism on $H$,
this induces the Gelfand triple
\[ V \subset H \subset V^*\]
such that in particular
\[_{V^*}\langle u, v\rangle_V = \langle  u, v\rangle _H \quad \forall\: u \in H, v\in V.\]
In addition we shall also assume that the inner product $\langle
.,.\rangle_V$ induces a closed quadratic form on $H$. This implies
the existence of a densely defined selfadjoint operator $\L :
H\supset D(\L) \to H$  on $H$ such that $V= D(\sqrt \L)$,
 $\la u, v\ra_V = \la u, \L v\ra_H$ for $u \in V, v \in D(\L)$  and such that the closure of $\L: V\supset D(\L) \to V^*$, still denoted by $\L$,
 defines an isometry. Moreover we assume that $L$ has discrete spectrum with corresponding eigenbasis  $(e_i)_{i\geq n}$, which will be  the case
 if e.g.\ the embedding $V\subset H$ is compact.

\noindent Let $(W(t))_{t\geq 0}$  be a cylindrical white noise on
some separable Hilbert space $(U, \langle .,.\rangle _U)$ defined on
some probability space $(\Omega,\mathbb P, \mathcal F)$ and let
$\mathcal F_t = \sigma(W_s, s\leq t)$ be the associated filtration.
For $X=H$ resp. $X=V$ we denote by $L_2(U,X)$   the class of
Hilbert-Schmidt mappings from $U$ to $X$, equipped with the
Hilbert-Schmidt norm $\|M\|_{L_2(U,X)}^2 = \sum\limits_{i\geq 1}
\langle M u_i, M u_i\rangle_X$, where $(u_i)_{i\geq 1}$ is some
orthonormal basis of $U$. Let
\[ A : V \to V^*, \quad \sigma: V \to L_2(U,V)\]
be measurable maps, then the existence and uniqueness result below
applies to $H$-valued It\^{o}-type stochastic differential equations
of the form
\begin{equation}\label{sde0}
 \left\{
\begin{array}{ll}
du(t)=Au(t)dt+\sigma(u(t))dW_t \\
u(0)=u_0 \in H.
\end{array}
\right.
\end{equation}
Below we shall work under the following  set of assumptions on the coefficients $A$ and $B$.

\begin{itemize}
 \item[(H1)] (Hemicontinuity) For all $u,\:v,\:x \in V$ the map
\[ \R \ni \lambda \to \subvstern\langle A(u+\lambda v ), x\rangle_V\]
is continuous.
\item [(H2)] (Weak monotonicity)
There exists $c_1\in \R$ such that for all $u$, $v\in V$
\begin{equation*}
2\ \subvstern \la Au-Av,u-v\ra_V +
\|\sigma(u)-\sigma(v)\|_{L_2(U,H)}^2 \leq c_1 \| u-v\|_H^2
\end{equation*}
\item [(H3)] (Lyapunov condition)  For  $n \in \mathbb N$, the operator $A$ maps   $H^n:=\mathop{\rm span}\{e_1, \dots, e_n\}\subset V$ into  $V$ and  there  exists a constant
 $c_2 \in \R$  such  that
 \begin{equation*}
2\   \la Au,u\ra_V + \|\sigma(u)\|_{L_2(U,V)}^2 \leq c_2(1+ \| u
\|_V^2 )\quad \forall u \in H^n, n\in \mathbb N.
\end{equation*}
\item [(H4)] (Boundedness) There exists a constant $c_3 \in \R$ such that
\[ \|A(u) \|\subvstern \leq c_3(1+\| u\|_V).\]
\end{itemize}

\begin{rem}{\normalfont Note that
(H3) replaces the standard coercivity assumption in
\cite{krylovrozovskii}
\begin{equation*}
2\   \subvstern \la Au,u\ra_V + \|\sigma(u)\|_{L_2(U,H)}^2 \leq c_2
\| u \|_H^2 -c_4 \|u\|_V^\alpha, \quad \forall v \in V \tag{A}
\end{equation*}
for some positive constant $c_4$ and  $\alpha >1$. Both conditions
(H3) and (A) yield the compactness of the Galerkin approximation in
$V$. Condition (A) is used indirectly by applying the finite
dimensional It\^{o} formula to the square of the $H$-norm. In our
case we use condition (H3)  directly by application of the finite
dimensional It\^{o} formula to the squared $V$-norm functional.}
\end{rem}

\noindent Basically, a solution to  \eqref{sde0} is a $V$-valued
process such that the equation holds in $V^*$ in integral form, c.f.
\cite{krylovrozovskii}. The following  precise definition is taken
from  \cite{Ro}.
  \begin{defi}\label{Definition}
A continuous $H$-valued $(\mathcal{F}_t)$-adapted process
$(u(t))_{t\in [0,T]}$ is called a solution of \eqref{sde0}, if for
its $dt\otimes \PP$-equivalence class $[u]$ we have $[u]\in
L^2([0,T]\times \Omega, dt\otimes \PP,V)$ and $\PP$-a.s.
\begin{equation*}
u(t)=u(0)+\int_0^tA(\bar{u}(s))\:ds+\int_0^t\sigma(\bar{u}(s))\:dW_s,\quad
t\in[0,T],
\end{equation*}
where $\bar{u}$ is any $V$-valued progressively measurable
$dt\otimes\PP$-version of $[u]$.
\end{defi}

Now we can state the main result of this section as follows.

\begin{theorem}\label{result}
Assume that conditions (H1)-(H4) hold, then for any initial data $u_0\in V$, there
exists a unique solution $u$ to \eqref{sde0} in the sense of
Definition \ref{Definition}. Moreover,
\begin{equation*}
\E\left(\sup\limits_{t\in[0,T]}\|u(t)\|_H^2\right)<\infty.
\end{equation*}
\end{theorem}

\begin{proof}
The proof follows the standard path of spectral Galerkin
approximation, the only difference towards \cite{krylovrozovskii,
Ro} is the compactness argument, c.f., lemma \ref{compactnesslemma}
below. To this aim let $(e_n)_{n\geq 1}$ be  an  orthonormal basis
in $H$ of eigenfunctions for the operator $\L:\: H \supset D(\L) \to
H$. Clearly $(e_n)_{n\geq 1}\subset V$ and the set
$\s\{e_n,\:\:n\geq 1\}$ is dense in $V$. Let $H_n\defeq
\s\{e_1,\cdots,e_n\}$ and define $P_n:\:V^{\ast}\rightarrow H_n$ by
$$
P_n y\defeq \sum\limits_{i=1}^n\:_{V^{\ast}}\la
y,e_i\ra_{V}e_i,\quad y\in V^{\ast}.
$$
\noindent Then we have $P_n|_H$ is just the orthogonal projection
onto $H_n$ in $H$. We shall define the family of $n$-dimensional
Brownian motions by setting
$$
W^n_t\defeq \sum\limits_{i=1}^n\la W_t,f_i\ra_{U}
f_i=\sum\limits_{i=1}^n B^i(t)f_i,
$$

\noindent where $(f_i)_{i\geq 1}$ is an orthonormal basis of the
Hilbert space $U$. We now consider the $n$-dimensional SDE
 \begin{equation}\label{n-sde}
 \left\{
\begin{array}{ll}
du^n(t)=P_n A
u^n(t)dt+P_n\sigma(u^n(t))dW^n_t \\
u^n(0,x)=P_n u_0(x) ,
\end{array}
\right\},
\end{equation}
which is identified with a corresponding SDE $dx(t)= b ^n (x(t)) dt
+ \sigma^n(x(t))d B^n_t$ in $\R^n$ via the isometric map $\R^n \to
H^n, x \to \sum _{i=1}^n x_i e_i$. By \cite[remark 4.1.2]{Ro}
conditions (H1) and (H2) imply the continuity of the fields $x\to
b^n(x)\in \R^n$ and $x\to \sigma^n(x)\in \R^{n\times n}$. Moreover,
assumption (H2) implies
\[
 2 \langle b^n (x) -b^n(y), x-y\ra_{\R^n} + |\sigma^n(x) - \sigma^n(y)|_{L_2(\R^n,\R^n)}^2 \leq c_1 |x-y|^2, \quad \forall x,y \in \R^n
\]
and, by the equivalence of norms on $\R^n$, (H3) gives the bound
\[  2\la b^n(x),x\ra + |\sigma^n(x)|_{L_2(\R^n,\R^n)} \leq c_5 (1+|x|^2),
\]
for some $c_5 \in \R$. Hence, equation \eqref{n-sde} is a weakly monotone and coercive equation in $\R^n$ which has a unique globally defined solution,
 cf.\ \cite[chapter 3]{Ro}.
\begin{lemma} \label{compactnesslemma}
Let $u^n$ be the solution to equation \eqref{n-sde}, then for any
$T>0$ we have
\begin{equation}
\sup\limits_{0\leq t\leq T}\E\|u^n(t)\|^2_V\leq (c_2 T+\E\bigl(
\|u_0\|^2_V\bigr))e^{c_2T}.
\end{equation}
\end{lemma}
\begin{proof}
Due to the definition of $P_n$ we may write
\begin{equation*}
\la u^n(t),e_i\ra=\la u^n(0),e_i\ra+\int_0^t\Big\la
\sum\limits_{k=1}^n\:_{V^{\ast}}\la
A(u^n(s)),e_k\ra_{V}e_k\:ds,e_i\Big\ra+\Big\la\int_0^tP_n\sigma(u^n(s))\:dW^n_s,e_i\Big\ra .
\end{equation*}

\noindent Hence, the It\^{o} formula in $\R^n$ yields
\begin{equation*}
\begin{split}
\|u^n(t)\|^2_V&=\|u_0^n\|_V^2+2\int_0^t\la P_n
A(u^n(s)),u^n(s)\ra_{V}\:ds+\int_0^t\|P_n\sigma(u^n(s))\|^2_{L_2(U_n,V)}\:ds\\&+
M^n(t),\quad t\in[0,T],
\end{split}
\end{equation*}
\noindent $\PP$-a.s, where $U_n := {\rm span}\,\{f_1, f_2, \cdots, f_n\} \subset U$ and
\begin{equation*}
M^n(t)\defeq 2\int_0^t\la u^n(s),P_n\sigma
(u^n(s))\:dW^n_s\ra_V,\quad t\in[0,T],
\end{equation*}
is a local martingale. We consider a sequence of $\mathcal{F}_t$-
stopping times $\tau_j$ with $\tau_j\uparrow +\infty$ as $j\to
+\infty$ and such that $\|u^n(t\wedge \tau_j)(\omega)\|_V$ is
bounded uniformly in $(t,\omega)\in [0,T]\times \Omega$,
$M^n(t\wedge \tau_j)$, $t\in [0,T]$ is a martingale for each $j\in
\NN$. Then we have
\begin{equation}
\begin{split}
\E\|u^n(t\wedge \tau_j)\|^2_V&=\E
\|u_0^n\|_V^2+2\int_0^{t}\E\ein_{[0,\tau_j]}\la P_n
A(u^n(s)),u^n(s)\ra_{V}\:ds\\&+\int_0^{t}\E\ein_{[0,\tau_j]}\|P_n\sigma(u^n(s))\|^2_{L_2(U_n,V)}\:ds.
\end{split} \label{normstep}
\end{equation}

\noindent Now using the definition of the operators  $A$ and $P_n$ we
can write
\begin{equation*}
\begin{split}
\la P_n A(u^n(s)),u^n(s)\ra_{V}&=\Big\la
\sum\limits_{i=1}^n\:_{V^{\ast}}\la
A(u^n(s)),e_i\ra_{V}e_i,u^n(s)\Big\ra_{V}\\
&=\sum\limits_{i=1}^n\:_{V^{\ast}}\la A(u^n(s)),e_i\ra_{V}\la
e_i,u^n(s)\ra_{V}.
\end{split}
\end{equation*}

\noindent Since   $u^n(t)\in H_n$ for $t\in [0,T]$ and $(e_n)_{n\geq
1}\subset V$ by assumption (H3) we can write
 $$
_{V^{\ast}}\la A(u^n(s)),e_i\ra_{V}=\la A(u^n(s)),e_i\ra_H,
$$
this yields
\begin{align*}
\la P_n A(u^n(s)),u^n(s)\ra_{V}&=\sum\limits_{i=1}^n\:\la
A(u^n(s)),e_i\ra_{H}\la e_i,u^n(s)\ra_{V}\\
&=\sum\limits_{i=1}^n\:\la A(u^n(s)),e_i\ra_{H}\lambda_i  \la
e_i,u^n (s)\ra_H
\end{align*}
where $\{\lambda_i \geq 0\}$ are the eigenvalues of the operator
$L$.

\noindent Therefore we have
\begin{equation*}
\la P_n A(u^n(s)),u^n(s)\ra_{V}=\la A(u^n(s)),u^n(s)\ra_{V}.
\end{equation*}

Hence, the operator $P_n$ may be dropped in the fist integral on the
right hand side term of \eqref{normstep} such that by the second
part of assumption (H3) \noindent
\begin{equation*}
\begin{split}
\E\|u^n(t\wedge \tau_j)\|^2_V&\leq\E \|u_0^n\|_V^2+c_2
\int_0^{t}(1+\E\|u^n\|_{V}^2) ds.
\end{split}
\end{equation*}

\noindent Hence letting $j\to +\infty$ and using Fatou's lemma we
obtain
\begin{equation*}
\E\|u^n(t)\|^2_V\leq \E
\|u_0^n\|_V^2+c_2\int_0^{t}(1+\E\|(u^n(s))\|_{V}^2)\:ds.
\end{equation*}

\noindent Now Gronwall's lemma yields
\begin{equation}\label{Gr}
\E\|u^n(t)\|^2_V\leq (c_2 T+\E
\|u_0^n\|_V^2)e^{c_2T}.
\end{equation}

\noindent For the estimate of $\E \|u_0^n\|_V^2$, we use the
definition of $P_n$ and write
\begin{equation*}
\begin{split}
\|u_0^n\|_V^2=\|P_n u_0\|^2_V&=\la
P_nu_0,P_nu_0\ra_V=\sum\limits_{i=1}^n\sum\limits_{j=1}^n
\:_{V^{\ast}}\la u_0,e_i\ra_{V}\la e_i,e_j\ra_{V}\:_{V^{\ast}}\la
u_0,e_i\ra_{V}\\
%&=\sum\limits_{i=1}^n\sum\limits_{j=1}^n \la u_0,e_i\ra_{H}\la \L e_i, e_j\ra_{H}\:\la u_0,e_i\ra_{H}\\
 &=\sum\limits_{i=1}^n \lambda_i \la  u_0,e_i\ra_{H}^2 \leq  \sum\limits_{i=1}^\infty \lambda_i \la  u_0,e_i\ra_{H}^2  =\|u_0\|_{V}^2.
\end{split}
\end{equation*}
 \end{proof}

From here all remaining arguments from \cite[chapter 4]{Ro} carry
over without change in order to complete the proof of the theorem. To make the paper self-contained  we briefly recall  the main steps. Let
\begin{equation*}
K\defeq L^2([0,T]\times \Omega,dt\otimes\PP,V)\quad \mbox{and} \quad
J\defeq L^2([0,T]\times \Omega,dt\otimes\PP,L_2(U,H)).
\end{equation*}

\noindent Due to the bound (H4) and the reflexivity of  $K$  we find a
subsequence $n_k\to +\infty$ such that $ u^{n_k}\to \bar{u}$ weakly in $K$ and weakly in $L^2([0,T]\times
\Omega,\:dt\otimes \PP, H)$,  $v^{n_k}\defeq A(u^{n_k})\to v$ weakly in $K^{\ast}$ and  $\theta^{n_k}\defeq P_{n_k}\sigma(u^{n_k})\to \theta$
weakly in $J$. Passing to the limit  in  \eqref{n-sde} one obtains  in $V^*$
\begin{equation}\label{u-formula}
u(t)\defeq u_0+\int_0^tv(s)\dd s+\int_0^t \theta(s)\dd W(s),\quad
t\in[0,T],
\end{equation}
and in particular $u=\bar{u}\:\:dt\otimes\PP$-a.e.  Now  the
 following It\^{o} formula for $\|u_t\|_H $ is crucial (c.f.\
\cite{krylovrozovskii}).

\begin{theorem} \label{itotheorem}
Let $u_0\in L^2(\Omega,\mathcal{F}_0,\PP,H)$ and $v\in
L^2([0,T]\times \Omega,dt\otimes\PP,V^{\ast})$, $\theta\in
L^2([0,T]\times \Omega,dt\otimes\PP,L_2(U,H))$, both progressively
measurable. Define the continuous $V^{\ast}$-valued process
\begin{equation*}
u(t)\defeq u_0+\int_0^t v(s)\:ds+\int_0^t\theta(s)\:dW_s,\quad
t\in[0,T].
\end{equation*}

\noindent If for its $dt\otimes\PP$-equivalence class $[u]$ we
have $[u]\in L^2([0,T]\times \Omega,dt\otimes\PP,V)$, then $u$
is an $H$-valued continuous $\mathcal{F}_t$-adapted process,
\begin{equation*}
\E\bigl(\sup\limits_{t\in[0,T]}\|u(t)\|_H^2\bigr)<\infty
\end{equation*}
and the following It\^{o}-formula holds for the square of its
$H$-norm $\PP$-a.s.
\begin{equation}\label{Ito-formula}
%\begin{split}
\|u(t)\|_H^2=\|u_0\|_H^2
%&
+2\int_0^t\Big(_{V^{\ast}}\la
v(s),\bar{u}(s)\ra_{V}+\|\theta(s)\|^2_{L_2(U,H)}\Big)\:ds
%\\&
+2\int_0^t\la
u(s),\theta(s)\:dW_s\ra
%\end{split}
\end{equation}
 for all $t\in[0,T]$, where $\bar u$ is any $V$-valued progressively measurable $dt\otimes\PP$-version of $[u]$.
\end{theorem}

In view of    \eqref{u-formula} this implies that   $u$ is continuous in $H$,  $\E (\sup\limits_{0\leq t\leq T}\|u(t)\|_H^2 )<+\infty$
and
\begin{equation}\label{product}
%\begin{split}
\E
%&
\left(e^{-{c_1}
t}\|u(t)\|_H^2\right)-\E\left(\|u_0\|^2_H\right)
%\\&
=\E\bigl(\int_0^te^{-{c_1}s}(2\:_{V^{\ast}}\la
v(s),\bar{u}(s)\ra_V+\|\theta(s)\|^2_{L_2(U,H)}-{c_1}\|u(s)\|^2_H)\dd
s\bigr).
%\end{split}
\end{equation}

\noindent An analogous formula holds true   for   $(u^{n_k}(t))_{t\geq 0}$. Hence, for
%it holds that
% \begin{equation*}
% \begin{split}
% \E&\left(e^{-{c_1}
% t}\|u^{n_k}(t)\|_H^2\right)-\E\left(\|u^{n_k}_0\|^2_H\right)
% %\\&=\E\left(\int_0^te^{-{c_1}s}(2\:_{V^{\ast}}\la
% %A(u^{n_k}(s)),u^{n_k}(s)\ra_V+\|P_{n_k}\sigma(u^{n_k})
% %\|^2_{L_2(U_n,H)}-{c_1}\|u^{n_k}(s)\|^2_H)\dd
% %s\right)\\
% \\&\leq \E\left(\int_0^te^{-{c_1}s}(2\:_{V^{\ast}}\la
% A(u^{n_k}(s)),u^{n_k}(s)\ra_V+\|\sigma(u^{n_k})\|^2_{L_2(U,H)}-{c_1}\|u^{n_k}(s)\|^2_H)\dd
% s\right).
% \end{split}
% \end{equation*}
%To see this apply the product rule to the precess $u^{n_k}$ stopped
%as soon as it reaches level $j$ and the let $j\to +\infty$.
 $\Phi\in K$, using (H2),
\begin{equation*}
\begin{split}
\E&\left(e^{-{c_1}
t}\|u^{n_k}(t)\|_H^2\right)-\E\left(\|u^{n_k}_0\|^2_H\right)
%\\&\leq
% \E\left(\int_0^te^{-{c_2}s}(2\:_{V^{\ast}}\la
%A(u^{n_k}(s)),u^{n_k}(s)\ra_V+\|\sigma(u^{n_k})\|^2_{L_2(U,H)}-{c_2}\|u^{n_k}(s)\|^2_H)\dd
%s\right)\\
\\&  \leq \E\left(\int_0^te^{-{c_2}s}\Big(2\:_{V^{\ast}}\la
A(\Phi(s)),u^{n_k}(s)\ra_V+2\:_{V^{\ast}}\la
A(u^{n_k}(s))-A(\Phi(s)),\Phi(s)\ra_V\right.\\
&- \|\sigma(\Phi(s))\|^2_{L_2(U,H)}+2\la
\sigma(u^{n_k}(s)),\sigma(\Phi(s))\ra_{L_2(U,H)}-2 {c_2}\la
u^{n_k}(s),\Phi(s)\ra_H+{c_2}\|\Phi(s)\|_H^2\Big)\dd s\Big).
\end{split}
\end{equation*}
Letting $k\to +\infty$ one concludes that for  every nonnegative $\psi\in
L^{\infty}([0,T],\R)$
\begin{equation}\label{est-apres monotonie}
\begin{split}
&\liminf\limits_{k\to+\infty}\E\left(\int_0^T
\psi(t)(e^{-{c_1}t}\|u^{n_k}(t)\|_H^2-\|u_0^{n_k}\|_H^2)\dd t\right)\\
&\leq \E\left(\int_0^T
\psi(t)\left(\int_0^te^{-{c_1}s}\Big(2\:_{V^{\ast}}\la
A(\Phi(s)),\bar{u}(s)\ra_V+2\:_{V^{\ast}}\la
v(s)-A(\Phi(s)),\Phi(s)\ra_V\right.\right.\\
&- \|\sigma(\Phi(s))\|^2_{L_2(U,H)}+2\la
\theta(s),\sigma(\Phi(s))\ra_{L_2(U,H)}-2 {c_1}\la
u(s),\Phi(s)\ra_H+{c_1}\|\Phi(s)\|_H^2\Big)\dd s\Big).
\end{split}
\end{equation}
\noindent On the other hand,  due to  the weak lower semicontinuity of the norm in $K$
% \begin{equation*}
% \begin{split}
% \E\left(\int_0^T\psi(t)\|\bar{u}(t)\|^2\dd
% t\right)&=\lim\limits_{k\to+\infty}\E\left(\int_0^T\la
% \psi(t)\bar{u}(t),u^{n_k}(t)\ra_H\dd t\right)\\
% &\leq \left(\E\int_0^T\psi(t)\|\bar{u}(t)\|^2\dd t\right)^{\frac 12}
% \lim\limits_{k\to+\infty}\left(\E\int_0^T
% \psi(t)\|u^{n_k}(t)\|^2_H\dd t\right)^{\frac 12}<+\infty.
% \end{split}
% \end{equation*}
%
% \noindent But $u=\bar{u}$ $\PP$-a.e., this yields
\begin{equation}
\E\left(\int_0^T\psi(t)\| {u}(t)\|_H^2\dd t\right)\leq
\liminf\limits_{k\to+\infty}\left(\E\int_0^T
\psi(t)\|u^{n_k}(t)\|^2_H\dd t\right).
\end{equation}

\noindent Combining this with \eqref{product}  and  \eqref{est-apres monotonie} one obtains
that
% \begin{equation*}
% \begin{split}
% &\E\left(\int_0^T
% \psi(t)(e^{-{c_2}t}\|u(t)\|_H^2-\|u_0\|^2)\dd t\right)\\
% &\leq \E\left(\int_0^T
% \psi(t)\left(\int_0^te^{-{c_2}s}\Big(2\:_{V^{\ast}}\la
% A(\Phi(s)),\bar{u}(s)\ra_V+2\:_{V^{\ast}}\la
% v(s)-A(\Phi(s)),\Phi(s)\ra_V\right.\right.\\
% &- \|\sigma(\Phi(s))\|^2_{L_2(U,H)}+2\la
% \theta(s),\sigma(\Phi(s))\ra_{L_2(U,H)}-2 {c_2}\la
% u(s),\Phi(s)\ra_H+{c_2}\|\Phi(s)\|_H^2\Big)\dd s\Big).
% \end{split}
% \end{equation*}
%
% \noindent Hence by making use of \eqref{product} we deduce that
\begin{equation} \label{final}
\begin{split}
\E\bigl(\int_0^T
\psi(t) & \int_0^te^{-{c_1}s} (2\:_{V^{\ast}}\la
v(s)-A(\Phi(s)),\bar{u}(s)-\Phi(s)\ra_V
\\&
+\|\sigma(\Phi(s))-\theta(s)\|^2_{L_2(U,H)}-{c_1}\|u(s)-\Phi(s)\|_H^2)\dd
s\dd t\bigr)\leq 0.
\end{split}
\end{equation}
\noindent Taking $\Phi=\bar{u}$ in \eqref{final} we obtain
$\theta=\sigma(\bar{u})$. By applying \eqref{final} to
$\Phi=\bar{u}-\varepsilon \tilde{\Phi} h$ for $\varepsilon>0$ and
$\tilde{\Phi}\in L^{\infty}([0,T]\times \Omega,\R)$, $h\in V$ and
dividing both sides by $\varepsilon$ and letting $\varepsilon\to 0$,
by (H2) and Lebesgue's theorem we get
\begin{equation*}
\E\left(\int_0^T \psi(t)\left(\int_0^te^{-{c_1}s}
\tilde{\Phi}(s)\Big(2\:_{V^{\ast}}\la v(s)-A(\bar{u}(s)),h\ra_V\dd
s\right)\dd t\right)\leq 0.
\end{equation*}

\noindent By the arbitrariness of $\psi$ and  $\tilde{\Phi}$ we
conclude that $v=A(\bar{u})$. \\

As for  the uniqueness consider two solutions  $u^{(1)}$ and $u^{(2)}$ of \eqref{sde0}
%in the sense of definition \ref{Definition},
with %
initial condition
$u_0^{(1)} \in V$ and  $u_0^{(2)}\in V$ respectively. Applying  theorem \ref{itotheorem}
to $u = u^{(1)} - u^{(2)}$ together with condition (H2) and Gronwall's lemma
\begin{equation}
 \E\|u^{(1)}(t)-u^{(2)}(t)\|_H ^2\leq \| u^{(1)}_0 -u^{(2)}_0\|_H^2e^{2c_1t}. \label{contractionest}
\end{equation}
This implies uniqueness of the solution for given initial state. Theorem \ref{result} is proved. \end{proof}

\subsection{\small \sc Generalized solutions for initial condition in $H$}

By means of    \eqref{contractionest} it is possible to construct a unique  generalized solution to \eqref{sde0} for initial condition in $u_0 \in H$. In particular this yields a unique Feller process on $H$ which extends the regular strong solutions of \eqref{sde0}.
\begin{proposition} \label{absttractgensol} Assume  (H1) - (H4) , then there exists a unique time homogeneous $H$-valued Markov process
$(\hat u_t^x, t \geq 0, x \in H)$ such that $t\to \hat u_t^x$ solves the SPDE \eqref{sde0} in the sense of definition $\eqref{Definition}$ whenever $x=u_0 \in V$.  Moerover, $(\hat u^x_t)$ induces a Feller semigroup on  $H$, i.e. the space  $C_b(H)$ of bounded continuous function on $H$ is invariant under the the operation $\varphi \to P_t\varphi$, where  $ P_t \varphi (x) = \E(\varphi (\hat u^x_t)) , x\in H$ for any $t \geq 0$.  \end{proposition}

\begin{proof}  For $x \in V \subset H$ define $t \to \hat u^x_t\in H$ as the unique solution to \eqref{sde0} with initial condition $u_0 =x$. For arbitrary $x\in H$, choose a sequence $(x_k)_k  $ in $  V$ such that $\| x_k - x \|_H \to 0$, then by \eqref{contractionest} the sequence of processes $(t\to \hat u _t^{x_k})_{k\in N}$
 is Cauchy in $ C([0,\infty); L^2(\Omega,H))$ with respect to the topology of locally uniform convergence and define $(t \to \hat u_t^x)$ as the unique limit. For $\varphi \in C_b(H)$ define $P_t\varphi (x)$ as above, then \eqref{contractionest} obviously yields
\begin{equation}  \E\|\hat u^{x}_t -\hat u^y_t\|_H^2  \leq e^{2c_1 t} \|x-y\|^2_H,\quad t\geq 0, \label{contracttwo}
\end{equation}
which implies that $P_t \varphi  \in C_b (H)$ for $\varphi \in
C_b(H)$. To prove that $(\hat u^x_t)^{x\in H}_{t \geq 0}$ is Markov,
by the monotone class argument  it suffices to show for all $x\in H$
\begin{equation}
 \E\bigl(\psi  (\hat u^x_t)\cdot  \varphi_1 (\hat u^x_{s_1})\cdots   \varphi_n (\hat u^x_{s_n})\bigr) =  \E\bigl( P_{t-s_n} \psi  (\hat u^x_{s_n})\cdot  \varphi_1 (\hat u^x_{s_1})\cdots   \varphi_n (\hat u^x_{s_n})\bigr), \label{markoveq}
\end{equation}
for any $0\leq s_1 \leq s_2 \cdots \leq s_n < t $ and $\varphi_1, \dots, \varphi_n, \psi \in C_b(H)\cap \Lip(H)$. By \eqref{contracttwo} we have
\[|P_{t-s}\varphi  (x) - P_{t-s}\varphi (y)| \leq e^{c_2(t-s)} \Lip (\varphi)\|x-y\|_H \quad \forall \, \varphi \in \Lip(H),\]
 hence will be enough to show \eqref{markoveq} for $x \in V$, where it follows by standard arguments from the uniqueness of solutions of \eqref{sde0},
 their adaptedness to the  filtration $\mathcal F_s$, $s\geq 0$, which for $s\leq t$  is independent of the sigma algebra of increments $\mathcal G_{s,t} =\sigma(W_\sigma-W_s; s\leq \sigma \leq t)$, c.f.,\ \cite[proposition 4.3.5]{Ro}.
%
%
% For simplicity let us  consider the case $n=1$ only, the case $n \geq 2$ being treated analogously. To this aim, for $v \in V$ and $s\geq 0$ let $ t\to \Phi_{s,t}^v$, $t\geq s$, denote the solution to the SDE
% \begin{equation*}\label{tsde}
%  \left\{
% \begin{array}{ll}
% d\Phi_{s,t} =A\Phi_{s,t} dt+\sigma(\Phi_{s,t})dW_t,  \quad t \geq s \\
% \Phi_{s,s}(0)=v
% \end{array}
% \right.,
% \end{equation*}
% which is adapted to the filtration generated by the Brownian increments $\mathcal G_{s,t} =\sigma(b^{i}\sigma-b^{i}_s; s\leq \sigma \leq t; i\in \mathbb N\}$ and therfore independent of   From the uniqueness of strong solutions to \eqref{sde0} in $L^2([0,T]\times \Omega, dt\otimes \mathbb P ; H)$ and the continuity of the maps $t\to u_t^v\in L^2(\Omega;H)$, resp. $t\to \Phi_{s,t}^v\in L^2(\Omega;H)$ it follows that
% \[  \hat u^{v}_{t} = \Phi_{s,t}(\hat u^v_{s}) \quad  \mathbb P\mbox{-a.s.}\]
This proves the existence of $(\hat u^x_t; t\geq 0; x\in H)$ as in the claim of the theorem. Trivially, uniqueness of $(\hat u_t^x)$ follows from  \eqref{contracttwo} which holds for any $H$-valued closure of solutions to equation \eqref{sde0}.
\end{proof}

\section{Application: Stochastic Curve Shortening Flow  in (1+1) Dimension}
% The motion of mean curvature is a well studied and rich object in
% geometric PDE theory for which a variety of methods have been
% developed c.f.\ e.g.\ \cite{MR1931534}. As a physical model it can
% be obtained  as sharp interface limit of the Allen-Cahn equation
% for the phase field of a binary alloy, describing the motion of the
% interface between the two phases. Stochastic mean curvature flow was
% derived heuristically  in e.g.,\ \cite{citeulike:2163102} as a
% refined model incoporating the influence of thermal noise. In the
% (d+1)-dimensional graph case the corresponding SPDE is of the form
% \begin{equation}
%  d u = \sqrt {1 + |\nabla u|^2}\, {\rm div}\bigl( \frac{\nabla
% u}{\sqrt {1 + |\nabla u|^2}} \bigr)\, dt+ B(u, \nabla u)  \delta W,  \label{grafspde}
% \end{equation}
% where $\delta$ stands for Stratonovich or It\^{o} differential,
% depending on the model. We not that equation \eqref{grafspde} is an
% example of degenerate stochastic parabolic problems.

%We refer to the recent work \cite{BDR} and reference therein for more examples of degenerate problems.

\subsection{\small \sc Strong solutions for $u_0 \in H^{1,2}([0,1])$ and smooth noise} Let us now show how we can treat the model rigorously in the  case $d=1$, which is also known as curve shortening flow,  using the results of the previous section. The simple but essential observation is that  for $d=1$  the drift operator in the SPDE \eqref{grafspde} above  may be written
\begin{equation}\label{drift} Au  = \frac{\partial^2 _x u}{1+(\partial _x u )^2} =\partial_x
( \arctan(\partial_x  u )),\end{equation} which fits into our
slightly modified \krylovrozovskii{}  framework.  To this aim let
\[  H_0^{1} ([0,1]) \subset L^2([0,1]) \subset H^{-1}([0,1]),\]
be the Gelfand triple, which is induced from the Dirichlet Laplacian $L = \Delta$ on $L^2([0,1])$.\\

For $u \in H^1_0([0,1])$,  let $Au\in H^{-1}([0,1])$ be defined by
\[ _{H^{-1}_0}\langle Au, v\rangle_{H^{1}_0} = - \int_{[0,1]} \arctan(\px u) \px v dx, \quad \forall v \in H^{1},\]
which is clearly hemicontinuous in the sense of condition (H1), due
to the continuity and uniform boundedness of
$\zeta\rightarrow\arctan\zeta$. Trivially $A$ is also bounded in the
sense of (H4) because
\begin{equation}
 \| A u \|_{H^{-1}([0,1])} = \sup_{v \in H^{1}_0([0,1]), \|v\|_{H^1_0}\leq 1}   \int_{[0,1]} \arctan(\px u) \px  v dx  \leq (\frac \pi 2)^{1/2}. \label{exboundedness}
\end{equation}

Moreover, by the monotonicity of $\arctan$
\begin{equation}
 _{H^{-1}} \la A u -A v , u -v \ra_{H^{1}} = - \int_{[0,1]} (\arctan(\px u) -\arctan(\px v))(\px u-\px v) dx \leq 0 \label{arctanmon}.
\end{equation}
The eigenvectors of $L = \Delta_0$ are $e_i= (x\to \sin(i 2\pi x)),
i\in \mathbb N$, hence  $Au = \px^2u /(1+(\px u)^2) \in
H^1_0([0,1])$ for any $u \in H^n =\mathop{\rm span}
\{e_1,\dots,e_n\}\subset {H^1_0([0,1])}$. Moreover,
\begin{equation} \langle  Au, u \rangle _{H^1_0}= -\int_{[0,1]}  \frac{\px^2u }{1+(\px u)^2} \px^2 u(x) dx \leq 0 \quad \forall u\in H^n. \label{energyestimate}
 \end{equation}

Let $(\phi_i)_{i\in \mathbb N }$  denote a sequence of linear independent Lipschitz functions on $[0,1]\times \R$ such that $\phi(0,y) = \phi(1,y) =0$
for all $y \in \R$ and such that the stronger regularity assumption \eqref{regcond} holds for the noise field, and let  furthermore $U$ denote the
Hilbert  space obtained from the closure of the  span of  $\{\phi_i, i \in \mathbb N\}$ with respect to the inner product
 $\langle \sum\limits_{i=1}^n{\lambda_i \phi_i}, \sum_{j=1}^m{\eta_j \phi_j}\rangle_U := \sum_{i=1}^{n \wedge m}  \lambda_i \eta_i$. \\

Define the diffusion operator $B: {H^1_0([0,1])} \to L(U,L^2([0,1])$ by
\[ B(u)[\phi](x) = \phi(x,u(x)) \in L^2([0,1])\]
Note that  $B(u)$ is in fact in $L_2(U, L^2([0,1]))$ since
\begin{align*}
 \|B(u)(\phi_i)\|_{L^2([0,1])}^2 &= \int_{[0,1]} \phi_i(x,u(x))^2 dx = \int_{[0,1]} |\phi_i(x,u(x))-\phi_i(0,u(0))|^2 dx \\
& \leq (\Lip(\phi_i))^2 \int_{[0,1]} (x^2+u^2(x))dx =    (\Lip(\phi_i))^2 (\frac 1 3 + \|u\|^2_{L^2([0,1])}), \end{align*}
such that
\begin{equation}
 \| B(u)\|_{L_2(U, L^2([0,1]))}^2 = \sum_{i}  \|B(u)(\phi_i)\|_{L^2([0,1])}^2 \leq  (\frac 1 3 + \|u\|^2_{L^2([0,1])}) \cdot \Lambda^2 \label{exhilbertschmidt}
\end{equation}

% \begin{align*}
%  \|B(u)(\psi)\|_{L^2([0,1])} &\leq(\frac 1 3 + \| u\|_{L^2([0,1])})^{1/2} \smallsum_{i} |\lambda_i | \, \Lip(\phi_i)  \\
% & \leq  (\frac 1 3 + \| u\|_{L^2([0,1])})^{1/2} ( \smallsum_i (\Lip(\phi_i)) ^2)^{1/2}  ( \smallsum_i \lambda_i ^2)^{1/2}  \\
% & \leq (\frac 1 3 + \| u\|_{L^2([0,1])})^{1/2}\, \Lambda \, \|\psi\|_{U} .  \end{align*}
Moreover,
\begin{align}
   \| B(u)-B(v) \|_{L_2(U, L^2([0,1]))}^2 & = \smallsum_{i} \| B(u)[\phi_i]- B(v)[\phi_i] \|_{L^2([0,1])}^2 \nonumber \\
& = \smallsum _i \int_{[0,1]} (\phi_i(x,u(x)) - \phi_i(x,v(x)))^2 dx \nonumber \\
&\leq \Lambda^2 \|u-v\|_{L^2([0,1])}^2
\label{h2cond}.
\end{align}
Similarly, $B(u)[\phi] \in {H^1_0([0,1])}$ for $u\in H^{1}_0([0,1])$, and by the chain rule for weakly differentiable functions,
\begin{align*}
 \|B(u)(\phi_i)\|_{H^1_0([0,1])}^2 &= \int_{[0,1]} (\px \phi_i(x,u(x)))^2 dx \\ \nonumber
&\leq (\Lip(\phi_i))^2\int_{[0,1]}(1+ |\px u(x)|^2) dx  =   (\Lip(\phi_i))^2 (1 + \|u\|^2_{H^1([0,1])}),  \end{align*}
which yields
\begin{equation}
 \| B(u)\|_{L_2(U, H^1([0,1]))}^2 = \smallsum_{i}  \|B(u)(\phi_i)\|_{H^1([0,1])}^2 \leq  ( 1 + \|u\|^2_{H^1([0,1])}) \cdot \Lambda^2  \label{growthcond}
 \end{equation}

In view of \eqref{exboundedness} -- \eqref{growthcond} we conclude
that the conditions (H1) -- (H4) are  satisfied in the given case
with constants $c_1=c_2=\Lambda^2$  and $c_3=\sqrt{ \pi/2}$. Hence,
by theorem \ref{result} we arrive at the following result.

\begin{theorem}\label{smcfresult}
Assume the regularity condition  \eqref{regcond} holds for the noise
field, then for any $T>0$ there is a  (up to $dt\otimes \mathbb
P$-equivalence in $[0,T]\times \Omega$) unique
$H^{1}_0([0,1])$-valued process $(u_t)_{t\in [0,T]}$ solving the
SPDE \eqref{themodel} in the sense of definition \ref{Definition}.
\end{theorem}

\subsection{\small \sc Generalized Markovian solution in $L^2([0,1])$ for non-smooth noise} Proposition
\ref{absttractgensol} readily yields generalized solutions for initial condition in $L^2([0,1])$ as follows.
\begin{proposition}  \label{markovsmcf1} Under condition $\eqref{regcond}$ there is a unique $L^2([0,1])$-valued
Markov process $(\hat u^x _t,  t \geq 0,x\in L^2([0,1]))$ such that
$t \to \hat u^x_t $ is a strong solution to the equation
\eqref{themodel} when $x=u_0 \in H^{1,2}_0([0,1])$. Moreover, $(\hat
u^x _t)_{t\geq 0}$ induces a Feller semigroup on $C_b(L^2([0,1]))$.
\end{proposition}

However, noticing that estimates \eqref{exhilbertschmidt} and \eqref{h2cond} remain true under the weaker regularity condition \eqref{kunitaregcond}, by  similar arguments as in the proof of proposition \ref{absttractgensol} we arrive at the following well-posedness result for the SPDE \eqref{themodel} under the Kunita-type regularity condition \eqref{kunitaregcond}.

\begin{proposition}  \label{markovsmcf2} Under condition $\eqref{kunitaregcond}$ there is a unique $L^2([0,1])$-valued
Markov process $(\hat u^x _t;  t \geq 0,x\in L^2([0,1]))$
such that for $x \in H^1_0([0,1])$, $ (u^x_t)_{t\geq 0} $
is the limit, in the sense of  locally uniform convergence on $C\bigl([0,\infty); L^2(\Omega,\mathbb P ; L^2([0,1]))\bigr)$, of
the strong solutions to the SPDE
\[d u^{(k)} =   \frac{ \partial _x^2  u^{(k)} }{ 1 + (\partial _x u^{(k)})^2} \, dt+ \sum_{i}^k \phi_i(.,u^{(k)}(.))\,db_t^i, \quad u_0^{k} =x.\]
Moreover,
\begin{equation}\label{monotonie}
E\| \hat u^x_t - \hat u^y_t\|_{L^2([0,1])}^2 \leq e^{\Lambda^2t}  \|
x - y\|_{L^2([0,1])}^2 \quad \forall x,y \in L^2([0,1]), t\geq 0.
\end{equation}
In particular, the induced semigroup, $P_t \varphi (x) =
\E(\varphi(\hat u^x_t))$ for measurable $\varphi: H \to \R$,   is
Feller.
\end{proposition}

\section{Ergodicity  for Stochastic Curve Shorting Flow with Additive Noise}
In this final section we show existence and uniqueness of an invariant measure for  the  generalized $L^2([0,1])$-valued solution $(\hat u^x_t; t\geq 0, x\in L^2([0,1]))$ obtained in proposition \ref{markovsmcf1}  for the SPDE \eqref{themodel} in the additive noise case, i.e.
when
\begin{equation} du = \frac{\partial_x ^2 u }{1+ (\partial_x u)^2} dt + Q dW_t,  \quad u(0)=u_0 \in H^{1,2}_0([0,1]), \label{constantdriftspde}\end{equation}
where $W$ is cylindrical white noise on some abstract Hilbert space $U$ and $Q \in L_2(U, H^{1,2}_0([0,1]))$.  As an example consider the case of $U=L^2([0,1])$ and $Q = (-\Delta)^{-\beta}$ for $\beta > 3/4$,  with $\Delta$ being the Dirichlet Laplacian on $[0,1]$.

\smallskip
Note also that for additive noise the condition (H2) is satisfied with $c_1 =0$.   As a consequence of \eqref{contractionest}, the Feller semigroup on $L^2$ induced from the generalized solutios $\hat u$ of  \eqref{themodel}  by $P_t \varphi (x) = \E(\varphi (\hat u^x_t))$ has the so-called $\textit {e}$-property \cite{KPS}, i.e.\ for all bounded Lipschitz continuous functions $\varphi : L^2 \mapsto \R$
\begin{equation}
|P_t\varphi(x) - P_t\varphi(y) | \leq \mbox{Lip}(\varphi) \|x - y \| \quad \forall x, y \in L^2. \label{eprop}
\end{equation}
%
% The proof of theorem \ref{ergodicthm} below is based on an application of a recent abstract ergodicty criterion for Feller semigroups satisfying the \textit{e}-property, due to   Komorowski, Peszat and Szarek \cite{KPS}.

\begin{theorem} \label{ergodicthm} \textit{Let $(P_t)_{t\geq 0}$ denote the Feller semigroup on $L^2([0,1])$ corresponding to the generalized solution to \eqref{constantdriftspde}, then $(P_t)$ is ergodic, i.\ e.\ there is a unique $(P_t)$-invariant probablity measure $\mu$ on $L^2([0,1])$. In particular,   $\lim_{t  \to \infty} \frac 1 t \int_0^t \langle P_t \varphi ,\nu\rangle = \langle \varphi ,\mu\rangle$ for any Borel probability measure $\nu \in \mathcal M_1 (L^2([0,1]))$ and any bounded continuous $\varphi :L^2([0,1]) \mapsto \R$.   }
\end{theorem}

Let  $ \qt(x, \cdot )  \defeq
\frac{1}{T}\int_0^T\mu_{\hat u_t}\:dt $,  where
$\mu_{\hat u_t}$ denotes the distribution at time $t$ of the generalized solution $\hat u^x_t$ to \eqref{constantdriftspde} with initial conditon $u_0=x\in L^2$.

%
% \noindent The following proposition shows tightness of the family of
% measure $\left\{\mu_T,\;T\geq1\right\}$.
%

\begin{proposition}\label{tight}
For any $x \in L^2$   the family of measures $
\left\{\qt(x,\cdot),\: T\geq 1\right\}$ is tight on $L^2([0,1])$.
\end{proposition}

\textit{Proof.} Assume first that $ \in H^{1,2}_0([0,1])$. In view of
\begin{equation*}
|\xi|-\alpha\leq\arctg\xi\cdot\xi\leq \beta+|\xi|,\quad \xi\in
\R,\quad \alpha,\:\beta>0
\end{equation*} it holds that
\noindent \begin{align}
 _{H^{-1}}\la Av,v\ra_{H^1}&=-\int_0^1 \arctg
(\partial _x v) \cdot \partial _x v \:dx\leq -\int_0^1|\partial _x v |\:dx+\alpha \nonumber \\
&\leq-c\|v\|_{W^{1,1}(0,1)}+\alpha
\label{recurrent}
\end{align}
for some $c>0$, by Poincar\'e  inequality.

\smallskip
Let  now   $t\to u_t$ be the solution to equation
\eqref{constantdriftspde} with regular initial condition $x=u_0 \in
H^{1,2}_0([0,1])$, then theorem \ref{itotheorem} holds. Hence  by
the It\^{o} formula for $\| u_t \|_{L^2([0,1]}^2$ and
\eqref{recurrent} we have
\begin{equation}
\begin{split}
\E\|u(t)\|^2&= \E\|u(0)\|^2+2\E\int_0^t\: _{V^{\ast}}\la
A(\ub(s)),\ub(s)\ra_{V}
ds+\E\int_0^t\|Q\|_{\mathcal{L}_{HS}(U,H)}^2\:ds\\
&\leq \E\|u(0)\|^2-c \:\E\int_0^t\|\ub(s)\|_{W^{1,1}(0,1)}+Dt
\end{split}
\end{equation}
where $D\defeq \alpha+\|Q\|_{\mathcal{L}_{HS}(U,H)}^2$.
% Thus we have
% $$
% \E\|u(t)\|^2+c\:\E\int_0^t\|\ub(s)\|_{W^{1,1}(0,1)}\:ds\leq
% \E\|u(0)\|^2+ Dt.
% $$
 In particular,
\begin{equation}
\E\left(\frac{1}{t}
\int_0^t\|\ub(s)\|_{W^{1,1}(0,1)}\:ds\right)\leq
\frac{1}{c}\Big(\E\|x\|^2+D\Big)\quad \forall  t\geq 1.
\label{chebyshevbound} \end{equation}
Since the functional $L^2 ([0,1]) \ni u \to \| u \|_{W^{1,1}(0,1)}
\in \R \cup \{ \infty\}$ has compact sublevel sets in $L^2([0,1])$,
the claim follows for regular initial condition $x=u_0 \in H^{1,2}_0([0,1])$.  \\

For the tightness of $\qt(x, .)$ with general $x \in L^2$ recall
(e.g.\ \cite[Remark on p. 49]{Par}) that it is sufficient (and
necessary) to find for arbitrary $\epsilon>0, \delta>0$   a finite
union of $\delta$-balls  $S_\delta = \bigcup_{\dotsk}
B_{\delta}(x_i) \subset  L^2$  such that
\[  \qt(x, S_\delta ) > 1-  \epsilon \quad \forall \,T>1. \]
To this aim choose $z \in B_{\delta \epsilon/4}(x)\cap
H^{1,2}_0(0,1)$ and a finite union of $\delta/2$-balls
$S_{\delta/2}= \bigcup\limits_{\dotsk} B_{\delta/2}(x_i) $ such that
$\qt(z,S_{\delta/2}) \geq 1 -\frac{\epsilon}{2}$. Let   $S_\delta =
\bigcup\limits_{\dotsk} B_{\delta }(x_i) $ and    choose a bounded
Lipschitz function $\varphi $ on $L^2$ with $\chi_{S_{\delta/2}}
\leq \varphi \leq \chi_{S_\delta}$ and $\mbox{Lip}(\varphi) \leq
\frac 2 \delta$. Hence,  using \eqref{eprop}, for all $T>1$
\begin{align*}  \qt(x, S_{\delta}) \geq \frac 1 T \int_0^T P_s \varphi (x) ds & \geq  \frac 1 T \int_0^T P_s \varphi (z) ds - \frac 2 \delta  \|x-z\| \\
 &
\geq \qt(z, S_{\frac \delta 2}) - \frac{2\|x -z\|}{\delta}>
1-\epsilon. \tag*{$\Box$}\end{align*}

\begin{lemma}\label{detlemm}
For $x\in L^2(0,1)$, let $(  v^ x(t))_{t\geq 0}$ the (generalized) solution of
\eqref{constantdriftspde} corresponding to $Q=0$. Then it holds
\[\lim\limits_{t\to +\infty}\|v^x(t)\|=0.
\]\end{lemma}
\begin{proof}
First, we consider the case where the initial data $v_0\in
C_0^{\infty}(0,1)$ (space of $C^{\infty}$-differentiable function
compactly supported in $[0,1]$. We set
$M\defeq\|v^\prime_0\|_{\infty}$ and define a function $h(t)$ with
\begin{itemize}
\item [(i)] $h$ is of class $C^{\infty}(\R)$ and satisfies
$$
h(t)=\arctan t \quad \mbox{for}\:\:|t|\leq M
$$
$$
|h(t)|\leq |t|,\quad t\in \R.
$$
\item[(ii)] $h^\prime$ is a bounded function on $\R$ satisfies
$\inf _{x\in \R}h^\prime(x)\geq \mu>0$ for a positive constant
$\mu$.
\item[(iii)] $h^{\prime\prime}$ is a bounded function on $\R$.
\end{itemize}

\noindent For $T>0$ fixed,  consider the equation
\begin{equation}
\label{h} \left\{
\begin{aligned}
dv(t) & = (h(v_x(t)))_x\,dt ,\\
v(0)& = v_0.\\
 \end{aligned}
\right.
\end{equation}
Following a similar argument as in \cite{taniguchi} and a maximum
principle for uniformly parabolic equation we can prove that the
classical solution $v$ of \eqref{h} satisfies
$$
\sup\limits_{0\leq t\leq T}\|v_x\|_{\infty}\leq M.
$$
\noindent Hence from the construction of $h$ we deduce that this
solution is also the solution of \eqref{constantdriftspde} with
$Q=0$ corresponding to the initial data $v_0\in C^\infty(0,1)$. Now
we remark that for the function $z\mapsto \arctan z$ we can write
$$
\arctan z=k(z)\cdot z\quad \mbox{for all $z\in \R$},
$$

for some positive decreasing function $k$ on $\R$. Therefore by
using the energy estimate for the function $v(t)$ we can write
\begin{equation}
\begin{split}
\frac 12 \frac{d}{dt}\|v(t)\|^2&=-\la \arctan v_x(t), v_x(t)\ra_{L^2(0,1)}\\
&\leq -\inf\limits_{z\in B(0,M)}k(z)\:\|v_x(t)\|^2\\ &\leq
-\inf\limits_{z\in B(0,M)}k(z)\:\|v(t)\|^2.
\end{split}
\end{equation}

Thus we obtain
$$
\|v(t)\|^2\leq e^{-2 t\inf\limits_{z\in B(0,M)}k(z)}\|v_0\|^2.
$$

This implies the statement of the lemma for regular initial datum
$v_0$. For general $v_0\in L^2(0,1)$ we proceed by approximation and
let $v_0^n$ a sequence of functions in $C^\infty _0(0,1)$ which
converges to $v_0$ in $L^2(0,1)$ for $n\to +\infty$. For $n\geq 0$
we denote by $v_n(t)$ the solution corresponding to the initial
condition $v_0^n$. By using the fact that $v_n(t)\to 0$ as $t\to 0$
and a triangle inequality argument we deduce the statement of the
lemma for general initial datum $v_0\in L^2(0,1)$.

\end{proof}

\begin{lemma}\label{stablem}
For $x\in L^2(0,1)$, let $(  \hat{v}^ x(t))_{t\geq 0}$ the
(generalized) solution of \eqref{constantdriftspde} corresponding to
$Q=0$. Then for every $x \in L^2, T>0$ and $\epsilon >0$, it holds
that \[ \PP(\|\hat u^x_T -\hat v^x_T\|< \epsilon ) >0.\]
\end{lemma}

\textit{Proof.}
 First we suppose that $x\in V$ and denote by
$(v^x_t)_{t\geq 0}$ the solution corresponding to
\eqref{constantdriftspde} with $Q=0$. We write
$$
z(t)= u(t)-v(t),\quad t\geq 0.
$$

Then the process $z(t)_{t\geq 0}$ solves the equation

\begin{equation*}
 \left\{
\begin{array}{lll}
dz(t)=(A u(t)-Av(t))dt+Q dW_t\\
z(0)=0.
\end{array}
\right.
\end{equation*}

We set
$$
z(t)=y(t)+QW_t.
$$
Then we have

$$
dy(t)=(Au(t)-Av(t))\:dt.
$$

Therefore,

\begin{equation*}
\begin{split}
\frac 1 2 \frac{d}{dt}\|y(t)\|^2&= \subvstern\la
Au(t)-Av(t),y(t)\ra_{V}\:dt\\
&=\subvstern\la Au(t)-Av(t),z(t)\ra_V\:dt-\subvstern\la
Au(t)-Av(t),QW_t\ra_{V}\\
&\leq 2\Big(\frac{\pi}{2}\Big)^{\frac 12}\|QW_t\|_V\leq
2\Big(\frac{\pi}{2}\Big)^{\frac 12} \|QW_t\|_V.
\end{split}
\end{equation*}

Where we used the monotonicity of $A$ and \eqref{exboundedness} to
obtain the estimate in the last line. Thus we deduce for $0\leq
t\leq T$

$$
\|y(t)\|\leq c\:T\sup\limits_{0\leq t\leq T}\|QW_t\|_V,
$$
for some positive constant $c$. We now use the splitting of
$z(\cdot)$ and the Poincar\'e inequality to obtain for $0\leq t\leq T$
\begin{equation}\label{zestimate}
\|z(t)\|\leq (c\: T+\frac 1 2 )\sup\limits_{0\leq t\leq T}\|QW_t\|_V.
\end{equation}

\noindent For the case where $x\in H$ we proceed by approximation
and use the uniform bound \eqref{exboundedness} to obtain the same
estimate as in \eqref{zestimate} for the process
$z(t)=\hat u^x(t)-\hat v^x(t)$, $x\in H$. Since $Q$ is a Hilbert-Schmidt operator from $U$ to $V$,  $(QW_t)_{t\geq 0}$ is a
  continuous Gaussian random process with values in $V$. Hence, for all $\delta>0$
$$
\PP\:\Big(\sup\limits_{0\leq t\leq T}\|QW_t\|_V<\delta\Big)>0.
$$

Now let $\ve>0$ and take $\delta>0$ such that $(cT +1/2)\delta <\ve$.
Then
\begin{equation*}
\PP\:\Big(\|z(t)\|<\ve\Big)>\PP\:\Big(\sup\limits_{0\leq t\leq
T}\|QW_t\|_V<\delta\Big)>0.\tag*{$\Box$}
\end{equation*}

\def\qt{Q^T}
\def\dotsk{{i=1, \dots, k}}

\begin{proposition} \label{lowerbd}  \textit{For  every $\delta >0$ and every $x\in L^2([0,1])$ it holds that
\[  \liminf_{T\to \infty } \qt(x, B_\delta(0)) >0.\]
}
\end{proposition}

\textit{Proof.} We proceed in three steps.  Let $\delta >0$ and $x \in L^2([0,1])$ be given.\\
\textbf{Step 1.}   For $R>0$ let $C_R = \{ u \in L^2| u \in
W^{1,1}_0(0,1) , \|u\|_{1,1} \leq R\} $, which is a compact subset
of $L^2([0,1])$. From \eqref{chebyshevbound}  and  Chebychev's
inequality we deduce
\[  \qt (0,L^2([0,1])\setminus C_R) \leq \frac {c}{R} \quad \forall\: T >1.
\]
Hence we may pick some $R>0$ such that $ \qt(0,C_R) >\frac 3 4  $ for all $T>1$.  From now we omit the subscript $R$, i.e. $C=C_R$. \\
\textbf{Step 2.}  Claim: There is some $\epsilon_1 >0$, a $\gamma_1 >0$  and a finite sequence $T_1, \cdots,  T_k $, $T_i >0$ such that
\[  \frac 1  k \sum _{i = 1, \dots, k} P_{T_i}(x,B_\delta(0)) >\gamma_1 \quad \forall\: x \in C_{\epsilon_1},
 \]
where $ C_{\epsilon_1} = \{ u \in L^2([0,1]) \,| \:\, d_{L^2}(u, C)
< \epsilon_1\}$ and $P_T(x,\cdot)$ the transition probability
corresponding to $(\hat u^x(t))_{t\geq 0}$ at time $T$. In fact, by
lemma \ref{detlemm} for each $x \in L^2([0,1])$ there exists a $T_x$
and a $r_x>0$ such that $\hat v^x_{T_x} \in B_{\delta/4}(0) $. For
$T>0$ and $\delta
>0$ let
\[ D(x, T, \delta)\defeq \PP\{   \| \hat v^x_T -\hat u^x_T\|_{L^2([0,1])} \leq  \delta\} ,\] which is strictly positive by lemma \ref{stablem}.  Hence it follows that $P_{T_x}(x,B_{\frac \delta 2}(0)) \geq D(x,T_x,\delta/4)=: \gamma_x >0$.
 Similarly as in the second part of proposition \ref{tight} we may use \eqref{eprop}  to deduce that for each $x \in L^2([0,1])$ there exists $r_x>0$
  such that $P_{T_x}(y,B_\delta(0)) > \gamma_x /2$ for all $y \in B_{r_x}(x)$.  Since $C$ is compact we may select a finite
  sequence $(x_i, r_i)$,  $\dotsk$, such that $C \subset \bigcup_\dotsk B(x_i, r_i)$. Setting $T_i := T_{x_i}$   the claim follows with  $\epsilon_1 := \min_\dotsk r_i  $  and $\gamma_1 := \min_\dotsk \gamma_i /2k$.

\textbf{Step 3:}  Choose $\rho >0$ such that
\[ \qt(x, C_{\epsilon_1}) > \frac 1 2 \quad \forall\: x \in B_{\rho}(0).\]
This is possible by a similar argument as in the second part proposition \ref{tight}. Finally,  by analogous reasons as in step 2, we may find  some $T_0>0$ and $\gamma_2>0$ such that $P_{T_0}(x,B_\rho(0))> \gamma_2$. \\

Hence,
\begin{align*}
\liminf_T  & \qt(x, B_\delta (0))  =   \liminf_{T} \frac 1 T  \int_0^T  P_{s} (x, B_\delta (0)) ds\\
 & = \liminf_{T} \frac 1 k \sum _{\dotsk} \frac 1 T \int_0^T  P_{s+T_i +T_0} (x, B_\delta (0)) ds \\
& = \liminf_{T}\frac 1 k \sum _{\dotsk}\frac 1 T \int_0^T \int_{L^2([0,1])}\int_{L_2([0,1])}  P_{T_i}(z,B_\delta(0)) P_{s}(y,dz) P_{T_0} (x, dy) ds \\
& \geq   \liminf_{T}  \frac 1 T \int_0^T \int_{B_\rho(0)}\int_{C_{\epsilon_1}}  \frac 1 k \sum _{\dotsk}  P_{T_i}(z,B_\delta(0)) P_{s}(y,dz) P_{T_0} (x, dy) ds\\
&  \geq  \gamma_1  \liminf_{T}  \frac 1 T \int_0^T \int_{B_\rho(0)}
P_{s}(y,C_{\epsilon_1}) P_{T_0} (x, dy) ds \intertext{which, by
Fatou's lemma is bounded from below by }
&  \geq  \gamma_1  \int_{B_\rho(0)}   \liminf_{T}  \frac 1 T \int_0^T   P_{s}(y,C_{\epsilon_1}) P_{T_0} (x, dy) ds\\
& =    \gamma_1  \int_{B_\rho(0)}   \liminf_{T} \qt(y,C_{\epsilon_1}) P_{T_0} (x, dy) ds \\
& > \frac 1 2 \gamma_1 P_{T_0}(x, B_\rho(0)) > \frac 1 2 \gamma_1
\gamma_2>0. \tag*{$\Box$}
\end{align*}
In view of \eqref{eprop} and proposition \ref{lowerbd}, Theorem
\ref{ergodicthm} is now a consequence of   \cite[Theorem 1]{KPS},
where
 $\mathcal T = L^2([0,1])$ according to proposition \ref{tight}.

\def\polhk#1{\setbox0=\hbox{#1}{\ooalign{\hidewidth
  \lower1.5ex\hbox{`}\hidewidth\crcr\unhbox0}}}
  \def\polhk#1{\setbox0=\hbox{#1}{\ooalign{\hidewidth
  \lower1.5ex\hbox{`}\hidewidth\crcr\unhbox0}}}

\end{document}